\documentclass[12pt]{article}

\usepackage{amsthm,amsmath,amssymb,latexsym}

\newtheorem{thm}{Theorem}[section]

\newtheorem{cor}[thm]{Corollary}

\newtheorem{lem}[thm]{Lemma}
\newtheorem{conj}[thm]{Conjecture}

\theoremstyle{remark}
\newtheorem*{rmk}{Remark}

\def\d{{\rm d}}

\def\pf{\noindent {\it Proof.} }

\numberwithin{equation}{section}
\renewcommand{\qed}{{\hfill\rule{4pt}{7pt}}\medskip}

\begin{document}
\begin{center}
{\Large\bf The Eulerian Distribution on Involutions is Indeed\\[10pt] Unimodal}
\end{center}
\vskip 2mm
\centerline{Victor J. W. Guo$^1$  and Jiang Zeng$^2$}

\begin{center}
Institut Camille Jordan,
Universit\'e Claude Bernard (Lyon I)\\
F-69622, Villeurbanne Cedex, France \\
{\tt $^1$guo@math.univ-lyon1.fr, http://math.univ-lyon1.fr/\textasciitilde guo}\\
{\tt $^2$zeng@math.univ-lyon1.fr, http://math.univ-lyon1.fr/\textasciitilde zeng}
\end{center}

\vskip 0.7cm {\small \noindent{\bf Abstract.}
Let $I_{n,k}$ (resp.~$J_{n,k}$) be the number of
involutions (resp.~fixed-point free involutions) of $\{1,\ldots,n\}$ with $k$ descents.
Motivated by Brenti's conjecture which states that
the sequence $I_{n,0}, I_{n,1},\ldots, I_{n,n-1}$ is log-concave, 
we prove that the two sequences $I_{n,k}$ and $J_{2n,k}$
are unimodal in $k$, for all $n$. Furthermore, we conjecture that there are nonnegative integers
$a_{n,k}$ such that 
\begin{align*}
\sum_{k=0}^{n-1}I_{n,k}t^k=\sum_{k=0}^{\lfloor (n-1)/2 \rfloor}a_{n,k}t^{k}(1+t)^{n-2k-1}.
\end{align*}
This statement is stronger than the unimodality of $I_{n,k}$ but is also interesting in 
its own right. }



\section{Introduction}
A sequence $a_0,a_1,\ldots,a_n$ of real numbers
is said to be {\it unimodal} if for some $0\leq j\leq n$ we have
$a_0\leq a_1\leq\cdots\leq a_j\geq a_{j+1}\geq\cdots \geq a_n$,
and is said to be {\it log-concave} if
$a_i^2\geq a_{i-1}a_{i+1}$ for all $1\leq i\leq n-1$. Clearly a log-concave
sequence
of \emph{positive} terms is unimodal. The reader is referred
to Stanley's survey~\cite{Stanley89} for the surprisingly
rich variety of methods to show that a sequence is log-concave or unimodal.
As noticed by Brenti~\cite{Brenti}, even though log-concave
and unimodality have one-line definitions, to prove the unimodality or log-concavity
of a sequence  can sometimes be a very difficult task requiring the use of intricate
combinatorial constructions or of refined mathematical tools.

Let $\mathfrak{S}_n$ be the set of all permutations of $[n]:=\{1,\ldots,n\}$.
We say that a permutation $\pi=a_1 a_2\cdots a_n\in\mathfrak{S}_n$ has a
{\it descent} at $i$ ($1\leq i\leq n-1$) if $a_i>a_{i+1}$.
The number of descents of $\pi$ is called its descent number and is
denoted by $\d(\pi)$. A statistic on $\mathfrak{S}_n$  is
said to be \emph{Eulerian},
if it is equidistributed with the descent number statistic. Recall that
the polynomial 
$$
A_n(t)=\sum_{\pi\in \mathfrak{S}_n} t^{1+\d(\pi)}=\sum_{k=1}^{n}A(n,k)t^k
$$
is called an {\it Eulerian polynomial}.
It is well-known that the {\it Eulerian numbers}
$A(n,k)$ ($1\leq k\leq n$) form a unimodal sequence, of which
several proofs  have been published: such as the analytical one by
showing that the polynomial $A_n(t)$ has only real zeros \cite[p. 294]{Co},
by induction based on the recurrence relation of $A(n,k)$ (see \cite{Kurtz}),
or by combinatorial techniques (see \cite{Gasharov, Stembridge}).

Let $\mathcal{I}_n$ be the set of all involutions
in $\mathfrak{S}_n$ and $\mathcal{J}_{n}$ the set of all fixed-point free involutions
in $\mathfrak{S}_{n}$. Define
\begin{align*}
I_n(t)&=\sum_{\pi\in \mathcal{I}_n} t^{\d(\pi)}=\sum_{k=0}^{n-1}I_{n,k}t^k,\\[5pt]
J_{n}(t)&=\sum_{\pi\in \mathcal{J}_{n}} t^{\d(\pi)}=\sum_{k=0}^{n-1}J_{n,k}t^k.
\end{align*}
The first values of these polynomials are given in Table~1.
\begin{table}[h]
\caption{The polynomials $I_{n}(t)$ and $J_{n}(t)$ for $n\leq 6$.\label{table:in6}}
\begin{center}
{\footnotesize
\begin{tabular}{|l|l|l|}
\hline
$n$  & $I_n(t)$ & $J_{n}(t)$  \\\hline
1    & 1   &  0 \\\hline
2    & $1+t$ & $t$ \\\hline
3    & $1+2t+t^2$ & 0  \\\hline
4    & $1+4t+4t^2+t^3$ & $t+t^2+t^3$  \\\hline
5    & $1+6t+12t^2+6t^3+t^4$ & 0  \\\hline
6    & $1+9t+28t^2+28t^3+9t^4+t^5$ &  $t+3t^2+7t^3+3t^4+t^5$ \\\hline
\end{tabular}
}
\end{center}
\end{table}

As one may notice from Table~\ref{table:in6} that the coefficients of
$I_n(t)$ and $J_{n}(t)$ are \emph{symmetric} and \emph{unimodal}
for $1\leq n\leq 6$. Actually, the
symmetries had been conjectured by Dumont and
were first proved by Strehl~\cite{Strehl}. Recently,
Brenti (see \cite{Dukes}) conjectured that the
coefficients of the polynomial $I_n(t)$ are \emph{log-concave} and
Dukes~\cite{Dukes} has obtained some partial results on the unimodality of
the coefficients of $I_n(t)$ and $J_{2n}(t)$. Note that, in contrast to
Eulerian polynomials $A_n(t)$, the polynomials $I_n(t)$ and
$J_{2n}(t)$ may have  \emph{non-real zeros}.

In this paper we will prove that for $n\geq 1$,
the two sequences $I_{n,0}, I_{n,1},\ldots, I_{n,n-1}$ and
$J_{2n,1}, J_{2n,2}, \ldots, J_{2n,2n-1}$ are
unimodal. Our starting point is
the known generating functions of polynomials $I_n(t)$ and $J_{n}(t)$:
\begin{align}
\sum_{n=0}^{\infty}I_n(t)\frac{u^n}{(1-t)^{n+1}}
&=\sum_{r=0}^{\infty}\frac{t^r}{(1-u)^{r+1}(1-u^2)^{r(r+1)/2}}, \label{eq:indt} \\[5pt]
\sum_{n=0}^{\infty}J_n(t)\frac{u^n}{(1-t)^{n+1}}
&=\sum_{r=0}^{\infty}\frac{t^r}{(1-u^2)^{r(r+1)/2}}, \label{eq:jndt}
\end{align}
which have been obtained by  D\'esarm\'enien and Foata~\cite{DF} and  Gessel
and Reutenauer~\cite{GR} using different methods.  We first derive linear
recurrence formulas for $I_{n,k}$ and $J_{2n,k}$ in the next section and
then prove the unimodality by induction in Section~3. We end this paper 
with further conjectures beyond the unimodality of the two sequences 
$I_{n,k}$ and $J_{2n,k}$.

\section{Linear  recurrence formulas for $I_{n,k}$ and $J_{2n,k}$}
Since the recurrence formula for the numbers
$I_{n,k}$ is a little more complicated
than $J_{2n,k}$, we shall first prove it for the latter.

\begin{thm}\label{lem:recj}
For $n\geq 2$ and $k\geq 0$, the numbers $J_{2n,k}$ satisfy
 the following recurrence formula:
\begin{align}
2nJ_{2n,k}
&=[k(k+1)+2n-2]J_{2n-2,k}+2[(k-1)(2n-k-1)+1]J_{2n-2,k-1}\nonumber\\[5pt]
&\quad +[(2n-k)(2n-k+1)+2n-2]J_{2n-2,k-2}. \label{eq:recj}
\end{align}
Here and in what follows $J_{2n,k}=0$ if $k<0$.
\end{thm}
\pf Equating the coefficients of $u^{2n}$ in \eqref{eq:jndt}, we obtain
\begin{equation}\label{fnt-jnt}
\frac{J_{2n}(t)}{(1-t)^{2n+1}}=\sum_{r=0}^{\infty}{r(r+1)/2+n-1\choose n}t^r.
\end{equation}
Since
$$
{r(r+1)/2+n-1\choose n}=\frac{r(r-1)/2+r+n-1}{n}{r(r+1)/2+n-2\choose n-1},
$$
it follows from \eqref{fnt-jnt} that
$$
\frac{J_{2n}(t)}{(1-t)^{2n+1}}
=\frac{t^2}{2n}\left(\frac{J_{2n-2}(t)}{(1-t)^{2n-1}}\right)''
+\frac{t}{n}\left(\frac{J_{2n-2}(t)}{(1-t)^{2n-1}}\right)'
+\frac{n-1}{n}\frac{J_{2n-2}(t)}{(1-t)^{2n-1}},
$$
or
\begin{align}
J_{2n}(t)&=\frac{t^2(1-t)^2}{2n}J_{2n-2}''(t)
+\left[\frac{(2n-1)t^2(1-t)}{n}+\frac{t(1-t)^2}{n}\right]J_{2n-2}'(t) \nonumber\\[5pt]
&\quad +\left[(2n-1)t^2+\frac{(2n-1)(1-t)t}{n}
+\frac{(n-1)(1-t)^2}{n}\right]J_{2n-2}(t) \nonumber\\[5pt]
&=\frac{t^4-2t^3+t^2}{2n}J_{2n-2}''(t)
+\left[\frac{(2-2n)t^3}{n}+\frac{(2n-3)t^2}{n}+\frac{t}{n}\right]J_{2n-2}'(t)
\nonumber\\[5pt]
&\quad +\left[(2n-2)t^2+\frac{t}{n}+\frac{n-1}{n}\right]J_{2n-2}(t).
\label{eq:j2n-der}
\end{align}
Equating the coefficients of $t^{n}$ in \eqref{eq:j2n-der} yields
\begin{align*}
J_{2n,k}&=\frac{(k-2)(k-3)}{2n}J_{2n-2,k-2}
-\frac{(k-1)(k-2)}{n}J_{2n-2,k-1}
+\frac{k(k-1)}{2n}J_{2n-2,k} \\[5pt]
&\quad +\frac{(2-2n)(k-2)}{n}J_{2n-2,k-2}+\frac{(2n-3)(k-1)}{n}J_{2n-2,k-1}
+\frac{k}{n}J_{2n-2,k} \\[5pt]
&\quad +(2n-2)J_{2n-2,k-2}+\frac{1}{n}J_{2n-2,k-1}
+\frac{n-1}{n}J_{2n-2,k}.
\end{align*}
After simplification, we obtain \eqref{eq:recj}. \qed

\begin{thm}For $n\geq 3$ and $k\geq 0$, the numbers $I_{n,k}$ satisfy
 the following recurrence formula:
\begin{align}
nI_{n,k}
&=(k+1)I_{n-1,k}+(n-k)I_{n-1,k-1}+[(k+1)^2+n-2]I_{n-2,k} \nonumber\\[5pt]
&\quad +[2k(n-k-1)-n+3]I_{n-2,k-1}+[(n-k)^2+n-2]I_{n-2,k-2}. \label{eq:reci}
\end{align}
Here and in what follows $I_{n,k}=0$ if $k<0$.
\end{thm}
\pf Extracting the coefficients of $u^{2n}$ in \eqref{eq:indt}, we obtain
\begin{equation}\label{eq:int-g}
\frac{I_{n}(t)}{(1-t)^{n+1}}
=\sum_{r=0}^{\infty}t^r\sum_{k=0}^{\lfloor n/2\rfloor}
{r(r+1)/2+k-1\choose k}{r+n-2k\choose n-2k}.
\end{equation}
Let
$$
T(n,k):={x+k-1\choose k}{y-2k\choose n-2k},
$$
and
$$
s(n):=\sum_{k=0}^{\lfloor n/2 \rfloor}T(n,k).
$$
Applying Zeilberger's algorithm, the Maple package
{\tt ZeilbergerRecurrence(T,n,k,s,0..n)} gives
\begin{equation}\label{eq:zeil}
(2x+y+n+1)s(n)+(y+1)s(n+1)-(n+2)s(n+2)= 0,
\end{equation}
i.e.,
$$
s(n)=\frac{y+1}{n}s(n-1)+\frac{2x+y+n-1}{n}s(n-2).
$$
When $x=r(r+1)/2$ and $y=r$, we get
\begin{equation}\label{g(n)}
s(n)=\frac{r+1}{n}s(n-1)+\frac{r(r-1)+3r+n-1}{n}s(n-2).
\end{equation}
Now, from \eqref{eq:int-g} and \eqref{g(n)} it follows that
\begin{align*}
\frac{nI_{n}(t)}{(1-t)^{n+1}}
&=t\left(\frac{I_{n-1}(t)}{(1-t)^{n}}\right)'
+\frac{I_{n-1}(t)}{(1-t)^{n}}
+t^2\left(\frac{I_{n-2}(t)}{(1-t)^{n-1}}\right)''
+3t\left(\frac{I_{n-2}(t)}{(1-t)^{n-1}}\right)'\\[5pt]
&\quad+(n-1)\frac{I_{n-2}(t)}{(1-t)^{n-1}},
\end{align*}
or
\begin{align}
nI_{n}(t)
&=(t-t^2)I_{n-1}'(t)+[1+(n-1)t]I_{n-1}(t)
+t^2(1-t)^2I_{n-2}''(t)\nonumber\\[5pt]
&\quad+t(1-t)[3+(2n-5)t]I_{n-2}'(t)
+(n-1)[1+t+(n-2)t^2]I_{n-2}(t).
\label{eq:ing-rec}
\end{align}
Comparing the coefficients of $t^k$ in both sides of \eqref{eq:ing-rec}, we obtain
\begin{align*}
nI_{n,k}
&=kI_{n-1,k}-(k-1)I_{n-1,k-1}
+I_{n-1,k}+(n-1)I_{n-1,k-1} \\[5pt]
&\quad +k(k-1)I_{n-2,k}-2(k-1)(k-2)I_{n-2,k-1}+(k-2)(k-3)I_{n-2,k-2} \\[5pt]
&\quad +3kI_{n-2,k}+(2n-8)(k-1)I_{n-2,k-1}-(2n-5)(k-2)I_{n-2,k-2} \\[5pt]
&\quad +(n-1)I_{n-2,k}+(n-1)I_{n-2,k-1}+(n-1)(n-2)I_{n-2,k-2},
\end{align*}
which, after simplification, equals the right-hand side of \eqref{eq:reci}.
\qed

\begin{rmk}
The recurrence formula \eqref{eq:zeil} can also be proved by hand as
follows. It is easy to see that the generating function of $s(n)$ is
\begin{equation}\label{eq:sn-gen}
\sum_{n=0}^{\infty}s(n)u^n=(1-u^2)^{-x}(1-u)^{-y-1}.
\end{equation}
Differentiating \eqref{eq:sn-gen} with respect to $u$ implies that
$$
\sum_{n=0}^{\infty}ns(n)u^{n-1}
=\left(\frac{2ux}{1-u^2}+\frac{y+1}{1-u}\right)(1-u^2)^{-x}(1-u)^{-y-1},
$$
consequently,
\begin{align}
(1-u^2)\sum_{n=0}^{\infty}ns(n)u^{n-1}
&=[(2x+y+1)u+y+1](1-u^2)^{-x}(1-u)^{-y-1}\nonumber \\[5pt]
&=[(2x+y+1)u+y+1]\sum_{n=0}^{\infty}s(n)u^{n}. \label{eq:sn-deri}
\end{align}
Comparing the coefficients of $u^{n+1}$ in both sides of \eqref{eq:sn-deri}, we obtain
$$
(n+2)s(n+2)-ns(n)=(2x+y+1)s(n)+(y+1)s(n+1),
$$
which is equivalent to \eqref{eq:zeil}.
\end{rmk}

Note that the right-hand side of \eqref{eq:recj} (resp.~\eqref{eq:reci}) is
invariant under the substitution $k\rightarrow 2n-k$ (resp.~$k\rightarrow n-1-k$),
provided that the sequence $I_{n-1,k}$ (resp.~$J_{2n-2,k}$) is symmetric.
Thus, by induction we derive immediately
the symmetry properties of $J_{2n,k}$ and $I_{n,k}$
(see \cite{DF, GR, Strehl}).

\begin{cor} For $n,k\in\mathbb{N}$, we have
$$
I_{n,k}=I_{n,n-1-k},\quad J_{2n,k}=J_{2n,2n-k}.
$$
\end{cor}

It would be interesting to find a combinatorial proof
of the recurrence formulas \eqref{eq:recj}
and \eqref{eq:reci}, since such a proof could hopefully lead to a combinatorial
proof of the unimodality of these two sequences.

\section{Unimodality of the sequences $I_{n,k}$ and $J_{2n,k}$}
The following observation is crucial in our inductive proof of
the unimodality of the sequences $I_{n,k}$ ($0\leq k\leq n-1$)
and $J_{2n,k}$ ($1\leq k\leq 2n-1$).
\begin{lem}\label{lem:sum-axi}
Let $x_0,x_1,\ldots,x_n$ and $a_0,a_1,\ldots,a_n$ be real numbers such that
$x_0\geq x_1\geq\cdots\geq x_n\geq 0$ and
$a_0+a_1+\cdots+a_k \geq 0$ for all $k=0,1,\ldots,n.$
Then
$$
\sum_{i=0}^{n}a_ix_i\geq 0.
$$
\end{lem}
Indeed, the above inequality follows from the identity:
$$
\sum_{i=0}^{n}a_ix_i=\sum_{k=0}^{n}(x_k-x_{k+1})(a_0+a_1+\cdots+a_k),
$$
where $x_{n+1}=0$.

\begin{thm}\label{thm:uni-jn}
The sequence $J_{2n,1},J_{2n,2},\ldots,J_{2n,2n-1}$ is unimodal.
\end{thm}
\pf By the symmetry of $J_{2n,k}$, it is enough to show that
$J_{2n,k}\geq J_{2n,k-1}$ for all $2\leq k\leq n$.
We proceed by induction on $n$. Clearly, the $n=2$ case is obvious.
Suppose the sequence $J_{2n-2,k}$ is unimodal in $k$.
By Theorem \ref{lem:recj}, one has
\begin{align}
2n(J_{2n,k}-J_{2n,k-1})
=A_0 J_{2n-2,k}+A_1J_{2n-2,k-1}+A_2J_{2n-2,k-2}+A_3 J_{2n-2,k-3},
\label{eq:jn-rec}
\end{align}
where
\begin{align*}
A_0&=k^2+k+2n-2, \quad
A_1=4nk-3k^2-6n+k+6,\\[5pt]
A_2&=3k^2+4n^2-8nk-5k+12n-4,\quad
A_3=3k-k^2+4nk-4n^2-8n.
\end{align*}
We have the following two cases:
\begin{itemize}
\item
If $2\leq k\leq n-1$, then
$$
J_{2n-2,k}\geq J_{2n-2,k-1}\geq J_{2n-2,k-2}\geq J_{2n-2,k-3}
$$
by the induction hypothesis, and clearly
\begin{align*}
&A_0\geq 0,\quad
A_0+A_1=2(k-1)(2n-k)+4 \geq 0, \\[5pt]
&A_0+A_1+A_2=(2n-k)^2-3k+8n \geq 0,\quad
A_0+A_1+A_2+A_3=0.
\end{align*}
Therefore, by Lemma~\ref{lem:sum-axi}, we have
$$
J_{2n,k}-J_{2n,k-1} \geq 0.
$$
\item If $k=n$, then
$$
J_{2n-2,n-1}\geq J_{2n-2,n}=J_{2n-2,n-2}\geq J_{2n-2,n-3}
$$
by symmetry and the induction hypothesis.
In this case, we have $A_1=(n-2)(n-3)\geq 0$ and thus the corresponding condition
of Lemma~\ref{lem:sum-axi} is satisfied.
Therefore, we have
$$
J_{2n,n}-J_{2n,n-1} \geq 0.
$$
\end{itemize}
This completes the proof. \qed

\begin{thm}\label{thm:uni-in}
The sequence $I_{n,0},I_{n,1},\ldots,I_{n,n-1}$ is unimodal.
\end{thm}
\pf By the symmetry of $I_{n,k}$, it suffices to show that
$I_{n,k}\geq I_{n,k-1}$ for all $1\leq k\leq (n-1)/2$. From Table \ref{table:in6},
it is clear that the sequences $I_{n,k}$ are unimodal in $k$ for $1\leq n\leq 6$.

Now suppose $n\geq 7$ and
the sequences $I_{n-1,k}$ and $I_{n-2,k}$ are unimodal in $k$.
Replacing $k$ by $k-1$ in \eqref{eq:reci}, we obtain
\begin{align}
nI_{n,k-1}
&=kI_{n-1,k-1}+(n-k+1)I_{n-1,k-2}+(k^2+n-2)I_{n-2,k-1} \nonumber\\[5pt]
&\quad +[2(k-1)(n-k)-n+3]I_{n-2,k-2}+[(n-k+1)^2+n-2]I_{n-2,k-3}. \label{eq:reci2}
\end{align}
Combining \eqref{eq:reci} and \eqref{eq:reci2} yields
\begin{align}
n(I_{n,k}-I_{n,k-1})
&=B_0I_{n-1,k}+B_1I_{n-1,k-1}+B_2I_{n-1,k-2} \nonumber\\[5pt]
&\quad +C_0I_{n-2,k}+C_1I_{n-2,k-1}+C_2I_{n-2,k-2}+C_3I_{n-2,k-3}, \label{eq:inkk-1}
\end{align}
where
\begin{align*}
B_0 &=k+1,\quad B_1=n-2k,\quad B_2=-(n-k+1),\\[5pt]
C_0 &=(k+1)^2+n-2,\quad C_1=2nk-3k^2-2k-2n+5, \\[5pt]
C_2 &=n^2-4nk+3k^2+4n-2k-5,\quad C_3=-(n-k+1)^2-n+2.
\end{align*}

Notice that $I_{n-1,k}\geq I_{n-1,k-1}\geq I_{n-1,k-2}$ for $1\leq k\leq (n-1)/2$.
By Lemma~\ref{lem:sum-axi}, we have
\begin{equation}\label{eq:b0ink}
B_0I_{n-1,k}+B_1I_{n-1,k-1}+B_2I_{n-1,k-2}\geq 0.
\end{equation}
It remains to show that
\begin{equation}\label{eq:c0ink}
C_0I_{n-2,k}+C_1I_{n-2,k-1}+C_2I_{n-2,k-2}+C_3I_{n-2,k-3}\geq 0,
\quad\forall\,1\leq k\leq (n-1)/2.
\end{equation}
We need to consider the following two cases:
\begin{itemize}
\item If $1\leq k\leq (n-2)/2$, then
$$
I_{n-2,k}\geq I_{n-2,k-1}\geq I_{n-2,k-2} \geq I_{n-2,k-3}
$$
by the induction hypothesis, and
\begin{align*}
&C_0=(k+1)^2+n-2\geq 0,\quad C_0+C_1=(2k-1)(n-k-1)+k+3\geq 0, \\[5pt]
&C_0+C_1+C_2=(n-k+1)^2+n-2\geq 0,\quad C_0+C_1+C_2+C_3=0.
\end{align*}
\item If $k=(n-1)/2$, then by symmetry and the induction hypothesis,
$$
I_{n-2,k-1} \geq I_{n-2,k}=I_{n-2,k-2} \geq I_{n-2,k-3}.
$$
In this case, we have $C_1=(n-3)(n-7)/4\geq 0$ for $n\geq 7$.
\end{itemize}
Therefore, by Lemma~\ref{lem:sum-axi} the inequality \eqref{eq:c0ink}
holds. It follows from \eqref{eq:inkk-1}--\eqref{eq:c0ink} that
$$
I_{n,k}-I_{n,k-1}\geq 0,\quad\forall\,1\leq k\leq (n-1)/2.
$$
This completes the proof. \qed
\section{Further remarks and open problems}
Since $I_{n,k}=I_{n,n-1-k}$, we can rewrite $I_n(t)$ as follows:
\begin{align*}
I_n(t)=\sum_{k=0}^{n-1}I_{n,k}t^k
=\begin{cases}\displaystyle\sum_{k=0}^{n/2-1}I_{n,k}t^k(1+t^{n-2k-1}),
&\text{if $n$ is even,}\\[15pt]
I_{n,(n-1)/2}t^{(n-1)/2}+\displaystyle\sum_{k=0}^{(n-3)/2}I_{n,k}t^k(1+t^{n-2k-1}),
&\text{if $n$ is odd.}
\end{cases}
\end{align*}
Applying the well-known formula
$$x^{n}+y^{n}
=\sum_{j=0}^{\lfloor n/2\rfloor}(-1)^{j}\frac{n}{n-j}{n-j\choose j}(xy)^j (x+y)^{n-2j},
$$
we obtain
\begin{align}
I_n(t)=\sum_{k=0}^{\lfloor (n-1)/2 \rfloor}a_{n,k}t^{k}(1+t)^{n-2k-1},
\label{eq:sym-int}
\end{align}
where
$$
a_{n,k}=
\begin{cases}\displaystyle\sum_{j=0}^{k}(-1)^{k-j}
\frac{n-2j-1}{n-k-j-1}{n-k-j-1\choose k-j}I_{n,j},&\text{if $2k+1<n$,}\\[15pt]
I_{n,k}+\displaystyle\sum_{j=0}^{k-1}(-1)^{k-j}\frac{n-2j-1}{n-k-j-1}
{n-k-j-1\choose k-j}I_{n,j},&\text{if $2k+1=n$.}
\end{cases}
$$

The first values of $a_{n,k}$ are given in Table~\ref{table:ank},
which seems to suggest the following conjecture.
\begin{conj}\label{conj:ank}
For $n\geq 1$ and $k\geq 0$, the coefficients $a_{n,k}$ are nonnegative integers.
\end{conj}
\begin{table}[!h]
\caption{Values of $a_{n,k}$ for $n\leq 16$ and $0\leq k\leq \lfloor (n-1)/2\rfloor$.\label{table:ank}}
{\scriptsize
\begin{center}
\begin{tabular}{|l|c|c|c|c|c|c|c|c|c|c|c|c|c|c|c|c|}
\hline
$k\setminus n$&1&2&3&4 & 5& 6& 7 & 8 & 9 & 10 & 11 & 12  & 13  & 14& 15 & 16  \\\hline
0      & 1 & 1 & 1 & 1 & 1& 1& 1 & 1 & 1 & 1  & 1  & 1   & 1   & 1 & 1  & 1\\\hline
1      &   &   & 0 & 1 & 2& 4& 6 & 9 & 12& 16 & 20 & 25  & 30  & 36& 42 & 49\\\hline
2      &   &   &   &   & 2& 6& 18& 39& 79& 141& 239& 379 & 579 & 849& 1211& 1680\\\hline
3      &   &   &   &   &  &  & 0 & 18& 78& 272& 722& 1716& 3626& 7160& 13206& 23263\\\hline
4      &   &   &   &   &  &  &   &   & 20& 124& 668& 2560& 8360& 23536& 59824& 139457\\\hline
5      &   &   &   &   &  &  &   &   &   &    & 32 & 700 & 4800& 24160& 95680& 325572\\\hline
6      &   &   &   &   &  &  &   &   &   &    &    &     & 440 &  5480& 44632& 257964\\\hline
7      &   &   &   &   &  &  &   &   &   &    &    &     &     &      & 2176& 44376\\\hline
\end{tabular}
\end{center}
}
\end{table}

Since the coefficients of $t^k(1+t)^{n-2k-1}$ are symmetric and unimodal with center
of symmetry at $(n-1)/2$, Conjecture~\ref{conj:ank}, is stronger than
the fact that the coefficients of $I_n(t)$ are symmetric and unimodal.
A more interesting question is to give a combinatorial interpretation of $a_{n,k}$.
Note that the Eulerian polynomials can be written as
$$
A_n(t)=\sum_{k=1}^{\lfloor (n+1)/2 \rfloor}c_{n,k}t^k(1+t)^{n-2k+1},
$$
where $c_{n,k}$ is the number of increasing binary trees on $[n]$ with
$k$ leaves and no vertices having left children only (see \cite{Branden,FS,Gasharov}).

We now proceed to derive a recurrence relation for $a_{n,k}$.
Set $x=x(t)=t/(1+t)^2$ and
$$
P_n(x)=\sum_{k=0}^{\lfloor (n-1)/2 \rfloor}a_{n,k}x^k.
$$
Then we can rewrite \eqref{eq:sym-int} as
\begin{equation}\label{eq:intpx}
I_n(t)=(1+t)^{n-1}P_{n}(x).
\end{equation}
Differentiating \eqref{eq:intpx} with respect to $t$ we get
\begin{align}
I_n'(t)&=(n-1)(1+t)^{n-2}P_{n}(x)+(1+t)^{n-1}P_{n}'(x)x'(t), \label{eq:intxt} \\[5pt]
I_n''(t)&=(n-1)(n-2)(1+t)^{n-3}P_{n}(x)+2(n-1)(1+t)^{n-2}P_{n}'(x)x'(t)\nonumber\\[5pt]
&\quad +(1+t)^{n-1}P_{n}''(x)(x'(t))^2+(1+t)^{n-1}P_{n}'(x)x''(t), \\[5pt]
x'(t)&=\frac{1-t}{(1+t)^3},\quad x''(t)=\frac{2t-4}{(1+t)^4}. \label{eq:intxt2}
\end{align}
Substituting \eqref{eq:intpx}--\eqref{eq:intxt2}
into \eqref{eq:ing-rec}, we obtain
\begin{align}
&\hskip -2mm n(1+t)^{n-1}P_n(x) \nonumber\\[5pt]
&=[1+(2n-2)t+t^2](1+t)^{n-3}P_{n-1}(x)+t(1-t)^2(1+t)^{n-5}P_{n-1}'(x)\nonumber\\[5pt]
&\quad+[-(t^2+14t+1)(1-t)^2+(1+6t-18t^2+6t^3+t^4)n
+4t^2n^2](1+t)^{n-5}P_{n-2}(x)\nonumber\\[5pt]
&\quad+[3t(t^2-4t+1)(1-t)^2+4t^2(1-t)^2n](1+t)^{n-7}P_{n-2}'(x)\nonumber\\[5pt]
&\quad+t^2(1-t)^4(1+t)^{n-9}P_{n-2}''(x). \label{eq:inxpnx}
\end{align}
Dividing the two sides of \eqref{eq:inxpnx} by $(1+t)^{n-1}$ and noticing that
$t/(1+t)^2=x$, after a little manipulation we get
\begin{align}
nP_n(x)&=[1+(2n-4)x]P_{n-1}(x)+(x-4x^2)P_{n-1}'(x) \nonumber\\[5pt]
&\quad +[(n-1)+(2n-8)x+4(n-3)(n-4)x^2]P_{n-2}(x) \nonumber\\[5pt]
&\quad +[3x+(4n-30)x^2+(72-16n)x^3]P_{n-2}'(x)+(x^2-8x^3+16x^4)P_{n-2}''(x).\nonumber
\end{align}
Extracting the coefficients of $x^k$  yields
\begin{align*}
na_{n,k}&=a_{n-1,k}+(2n-4)a_{n-1,k-1}+ka_{n-1,k}-4(k-1)a_{n-1,k-1}\\[5pt]
&\quad+(n-1)a_{n-2,k}+(2n-8)a_{n-2,k-1}+4(n-3)(n-4)a_{n-2,k-2}\\[5pt]
&\quad+3ka_{n-2,k}+(4n-30)(k-1)a_{n-2,k-1}+(72-16n)(k-2)a_{n-2,k-2} \\[5pt]
&\quad+k(k-1)a_{n-2,k}-8(k-1)(k-2)a_{n-2,k-1}+16(k-2)(k-3)a_{n-2,k-2}.
\end{align*}
After simplification, we obtain the following recurrence formula for $a_{n,k}$.

\begin{thm}\label{thm:rec-ank}
For $n\geq 3$ and $k\geq 0$, there holds
\begin{align}
na_{n,k}
&=(k+1)a_{n-1,k}+(2n-4k)a_{n-1,k-1}+[k(k+2)+n-1]a_{n-2,k}  \nonumber\\[5pt]
&\quad+[(k-1)(4n-8k-14)+2n-8]a_{n-2,k-1}+4(n-2k)(n-2k+1)a_{n-2,k-2}, \label{eq:rec-ank}
\end{align}
where $a_{n,k}=0$ if $k<0$ or $k>(n-1)/2$.
\end{thm}

Note that, if $n\geq 2k+3$, then $$(k-1)(4n-8k-14)+2n-8>0
\quad\text{for any $k\geq 1$,}$$ and so are the other coefficients
in \eqref{eq:rec-ank}. Therefore, 
Conjecture~\ref{conj:ank} would be proved if one can show that
$a_{2n+1,n}\geq 0$ and $a_{2n+2,n}\geq 0$.

{}Finally, from \eqref{eq:sym-int} it is easy to see that
\begin{align*}
a_{2n+1,n}&=(-1)^nI_{2n+1}(-1)=\sum_{k=0}^{2n}(-1)^{n-k} I_{2n+1,k},\\[5pt]
a_{2n+2,n}&=(-1)^{n}I_{2n+2}'(-1)=\sum_{k=1}^{2n+1}(-1)^{n+1-k} kI_{2n+2,k}.
\end{align*}
Thus, Conjecture~\ref{conj:ank} is equivalent to the {\it nonnegativity}
of the above two alternating sums.

Since $J_{2n,k}=J_{2n,2n-k}$, in the same manner as $I_{n}(t)$ we obtain
\begin{align*}
J_{2n}(t)=\sum_{k=1}^{n}b_{2n,k}t^{k}(1+t)^{2n-2k},
\end{align*}
where
$$
b_{2n,k}=
\begin{cases}\displaystyle\sum_{j=1}^{k}(-1)^{k-j}
\frac{2n-2j}{2n-k-j}{2n-k-j\choose k-j}J_{2n,j},&\text{if $k<n$,}\\[15pt]
J_{2n,k}+\displaystyle\sum_{j=1}^{k-1}(-1)^{k-j}\frac{2n-2j}{2n-k-j}
{2n-k-j\choose k-j}J_{2n,j},&\text{if $k=n$.}
\end{cases}
$$
Now, it follows from \eqref{fnt-jnt} that
\begin{align*}
J_{2n,k}=\sum_{i=0}^{k}(-1)^{k-i}{2n+1\choose k-i}{i(i+1)/2+n-1\choose i(i+1)/2-1}
\end{align*}
is a polynomial in $n$ of degree $d:=k(k+1)/2-1$ with leading coefficient
$1/d!$, and so is $b_{2n,k}$. Thus, we have
$\lim_{n\rightarrow+\infty} b_{2n,k}=+\infty$ for any fixed $k>1$.

The first values of $b_{2n,k}$ are given in Table~\ref{table:bnk},
which seems to suggest
\begin{conj}\label{conj:bnk}
For $n\geq 9$ and $k\geq 1$, the coefficients $b_{2n,k}$
are nonnegative integers.
\end{conj}
\begin{table}[h]
\caption{Values of $b_{2n,k}$ for $2n\leq 24$ and
$1\leq k\leq n$.\label{table:bnk}}
{\scriptsize
\begin{center}
\begin{tabular}{|l|c|c|c|c|c|c|c|c|c|c|c|c|}
\hline
$k\setminus 2n$
  &2& 4& 6& 8&   10& 12&  14&   16&    18&     20&      22&       24  \\\hline
1 &1& 1& 1& 1&  1&  1&   1&    1&     1&      1&       1&        1  \\\hline
2 & &$-1$&$-1$ &0&   2&  5&   9&   14&    20&     27&      35&       44  \\\hline
3 & &  & 3&12& 36& 91& 201&  399&   728&   1242&    2007&     3102  \\\hline
4 & &  & &$-7$&$-10$ &91&652&2593& 7902& 20401&   46852&    98494  \\\hline
5 & &  &  &  & 25&219&1710&10532& 50165& 194139&  639968&  1861215  \\\hline
6 & &  &  &  & &$-65$& 249&11319&122571& 841038& 4377636& 18747924  \\\hline
7 & &  &  &  &   &   & 283& 6586&135545&1737505&15219292&101116704  \\\hline
8 & &  &  &  &   &   &   &$-583$& 33188&1372734&24412940&277963127  \\\hline
9 & &  &  &  &   &   &    &     &  4417& 379029&16488999&367507439  \\\hline
10& &  &  &  &   &   &    &     &      &   1791& 3350211&203698690  \\\hline
11& &  &  &  &   &   &    &     &      &       &  133107& 36903128  \\\hline
12& &  &  &  &   &   &    &     &      &       &        &   761785  \\\hline
\end{tabular}
\end{center}
}
\end{table}

Similarly to the proof of Theorem~\ref{thm:rec-ank}, we can prove the following
result.
\begin{thm} \label{thm:rec-bnk}
For $n\geq 2$ and $k\geq 1$, there holds
\begin{align*}
2n b_{2n,k}&=[k(k+1)+2n-2]b_{2n-2,k}+[2+2(k-1)(4n-4k-3)]b_{2n-2,k-1}\\[5pt]
&\quad+8(n-k+1)(2n-2k+1)b_{2n-2,k-2}.
\end{align*}
where $b_{2n,k}=0$ if $k<1$ or $k>n$.
\end{thm}

Theorem \ref{thm:rec-bnk} allows us to reduce the verification of Conjecture~\ref{conj:bnk}
to the boundary case $b_{2n,n}\geq 0$ for $n\geq 9$.

\vskip 5mm
\noindent{\bf Acknowledgment.}
The second author was supported by EC's IHRP Programme, within Research Training
Network ``Algebraic Combinatorics in Europe," grant HPRN-CT-2001-00272.

\renewcommand{\baselinestretch}{1}

\end{document}